\documentclass[10pt,twoside]{article}
\usepackage{amsmath,amsthm}
\usepackage[latin1]{inputenc}
\usepackage{picinpar}
\usepackage{longtable}
\usepackage{amsrefs}
\usepackage[dvips]{graphicx}
\usepackage[mathscr]{eucal}
\usepackage{calc}
\usepackage{ifthen}
\usepackage{enumerate}
\usepackage{amssymb,latexsym}
\usepackage{amsthm}
\usepackage[dvips]{graphics}
\usepackage{color}
\usepackage{tabularx}
\usepackage{makeidx}
\theoremstyle{definition}

 \newtheorem{teor}{Theorem}

\newtheorem{Nota}{Remark}

\begin{document}
\thispagestyle{plain}
\par\bigskip
\begin{centering}

\textbf{On sums of three squares}

\end{centering}\par\bigskip
 \begin{centering}
 \footnotesize{AGUSTIN MORENO CA\~NADAS}\par\bigskip
 \end{centering}

$\centerline{\textit{\small{Departamento de Matem\'aticas, Universidad Nacional de Colombia}}}$
\centerline{\textit{\small{Bogot\'a-Colombia}}} $\centerline{\textit{\small{amorenoca@unal.edu.co}}}$

\bigskip

\bigskip \small{We prove that a positive integer not of the form $4^{k}(8m+7)$, $k,m\in\mathbb{N}$ can be expressible as a sum of three or fewer
squares by using some results of Kane and Sun on mixed sums of squares and triangular numbers.}
\par\bigskip
\small{\textit{Keywords} : quadratic form, representation, square number, triangular number.}

\bigskip \small{Mathematics Subject Classification 2000 : 05A17 ; 11D45; 11D85; 11E25; 11P83.}

\bigskip

\textbf{1. Introduction}\par\bigskip

One of the most investigated topics in additive number theory is the representation of integers by sums of squares
and more generally, by quadratic forms. For example the classical problem of finding formulas for the number of
ways of expressing an integer as the sum of $s$ squares. One can also ask for every number to be expressible as
the sum of as few as possible square numbers. For instance there is Gauss's famous 1796-07-10 diary
entry\par\bigskip

\begin{centering}
E$\Upsilon$PHKA!\qquad $\mathrm{num}=\Delta+\Delta+\Delta$, \par\bigskip

\end{centering}

that is, Gauss proved that every natural number is the sum of three or fewer triangular numbers. This statement is
equivalent to the statement that every number of the form $8m+3$ is a sum of three odd squares. Actually Gauss's
theorem implies the Lagrange's theorem (1772), that every natural number is a sum of four or fewer square numbers
[4,10]. \par\bigskip

Legendre proved in 1798 that the set of positive integers that are not sums of three or fewer squares
$=\{n\in\mathbb{N}\backslash\{0\}\mid
n=4^{s}(8m+7),\hspace{0.1cm}\mathrm{for\hspace{0.1cm}some}\hspace{0.1cm}m,s\in\mathbb{N}\}$. Shortly afterwards,
in 1801, Gauss going way beyond Legendre, actually obtained a formula for the number of primitive representations
of an integer as a sum of three squares. According to Ewell [4,8] and others authors no simple proof of this
theorem has been found up to date.\par\bigskip

At the present time, we know that Lagrange's theorem is a particular case of the fifteen theorem of Conway and
Schneeberger, which states that if a positive integer-matrix quadratic form represents each of 1, 2, 3, 5, 6, 7,
10, 14, 15, then it represents all positive integers [1-3]. Bhargava gave a simple proof of this theorem [1], and
Kane proved a similar condition for sums of triangular numbers [13].\par\bigskip

The following more general theorem (290-theorem) was proved by Bhargava and Hanke [2,3].\par\bigskip If a
positive-definite integral quadratic form represents each of \par\smallskip
1,2,3,5,6,7,10,13,14,15,17,19,21,22,23,26,29,30,31,34,35,37,42,58,93,110,145,203,290,\par\smallskip then it
represents all positive integers. \par\bigskip

On partitions into square numbers, Jacobi by the use of elliptic and theta functions proved that the number of
representations of a positive integer $n$ as the sum of four squares is given by $8[2+(-1)^{n}]\sigma_{0}$, where
$\sigma_{0}$ denotes the sum of the odd divisors of $n$. Lehmer denoted $P_{k}(n)$ the number of partitions of a
natural number $n$ into $k$ integral squares $\geq0$, and solved almost completely the equation $P_{k}(n)=1$ [12].
Lehmer claimed that the general problem of finding a formula for $P_{k}(n)$ was a problem of great complexity. The
case $k=3$ was studied by Grosswald, A. Calloway, and J. Calloway in [6], and Grosswald solved (essentially) the
problem, giving the number of partitions of an arbitrary integer $n$ into $k$ squares (taking into account that,
he didn't distinguish between partitions that contains zeros and those that do not) [7]. \par\bigskip

In this paper we shall give a solution to the following problem proposed by Guy in [10] :\par\bigskip

What theorems are there, stating that all numbers of a suitable shape are expressible as the sum of three squares
of numbers of a given shape? \par\bigskip

In order to obtain a proof of the difficult part of (it is easy to verify the only if part) Legendre-Gauss
theorem we will use the solution to the problem described above and some new results of Kane and Sun on almost
universal mixed sums of squares and triangular numbers.\par\bigskip

\textbf{1. On sums of squares and triangular numbers}\par\bigskip

In this section we describe some recent results concerning representations of numbers by sums of triangular and
square numbers.\par\bigskip

We let $t_{k}=\frac{k(k+1)}{2}$, $s_{k}=k^{2}$ denote the triangular and square $k$-th numbers
respectively.\par\bigskip

The following theorem proved by Lebesgue and R\'{e}alis in [15] was reproved by Farkas in [5], via the theory of
theta functions,

\begin{teor}\label{Farkas}

Every positive integer can be written as the sum of two squares plus one triangular number and every positive
integer can be written as the sum of two triangular numbers plus one square.

\end{teor}

In [4] Ewell proved the following theorem

\begin{teor}\label{Ewell}

For each $n\in\mathbb{N}$, $t_{2}(n)=d_{1}(4n+1)-d_{3}(4n+1)$.

\end{teor}

Where $t_{2}(n)$ is the number of representations of $n$ by sums of 2 triangular numbers and $d_{i}(n)$ is the
number of positive divisors of $n$ congruent to $i\hspace{0.1cm}\mathrm{mod}\hspace{0.1cm}4$.\par\bigskip In
[9,16], Guo, Pan, and Sun showed the following theorem
\begin{teor}\label{Sun}
\begin{enumerate} [(a)]
\item Any natural number is a sum of an even square and two triangular numbers, and each positive integer is a sum
of a triangular number plus $s_{x}+s_{y}$ for some $x,y\in\mathbb{Z}$ with $x\not\equiv
y\hspace{0.1cm}\mathrm{mod}\hspace{0.1cm}2$ or $x=y>0$.

\item Let $a,b,c$ be positive integers with $a\leq b$. Every $n\in\mathbb{N}$ can be written as $as_{x}+bs_{y}+ct_{z}$ with $x,y,z\in\mathbb{Z}$ if and only if $(a,b,c)$ is among the following vectors :\par\bigskip

$(1,1,1)$, $(1,1,2)$, $(1,2,1)$, $(1,2,2)$, $(1,2,4)$,\par\smallskip $(1,3,1)$, $(1,4,1)$, $(1,4,2)$, $(1,8,1)$,
$(2,2,1)$.\par\bigskip

\item Let $a,b,c$ be positive integers with $b\geq c$. Every $n\in \mathbb{N}$ can be written as $as_{x}+bt_{y}+ct_{z}$
with $x,y,z\in\mathbb{Z}$ if and only if $(a,b,c)$ is among the following vectors :\par\bigskip

$(1,1,1)$, $(1,2,1)$, $(1,2,2)$, $(1,3,1)$, $(1,4,1)$, $(1,4,2)$, $(1,5,2)$,\par\smallskip $(1,6,1)$, $(1,8,1)$,
$(2,1,1)$, $(2,2,1)$, $(2,4,1)$, $(3,2,1)$, $(4,1,1)$, $(4,2,1)$.

\end{enumerate}
\end{teor}

In [13] Kane gave the following generalization of Gauss's Eureka theorem,

\begin{teor}\label{Kane}

Fix the sequence $b_{1}, b_{2},\dots,b_{k}$. Then

\begin{enumerate}[(a)]

\item The sum of triangular numbers\par\bigskip
\begin{centering}

$f(x)=f_{b}(x)=\underset{i=1}{\overset{k}{\sum}}b_{i}t_{x_{i}}$\par\bigskip
\end{centering}

represents every positive integer if and only if $f_{b}$ represents the integers $1,2,4,5,$ and 8.

\item The corresponding diagonal quadratic form $Q(x)=\underset{i=1}{\overset{k}{\sum}}b_{i}s_{x_{i}}$ with
$x_{i}$ all odd represents every integer of the form\par\bigskip

\begin{centering}
$8n+\underset{i=1}{\overset{k}{\sum}}b_{i}$\par\bigskip

\end{centering}

if and only if it represents $8+\underset{i=1}{\overset{k}{\sum}}b_{i}$,
$16+\underset{i=1}{\overset{k}{\sum}}b_{i}$, $32+\underset{i=1}{\overset{k}{\sum}}b_{i}$,
$40+\underset{i=1}{\overset{k}{\sum}}b_{i}$, and $64+\underset{i=1}{\overset{k}{\sum}}b_{i}$.

\end{enumerate}

\end{teor}

Kane and Sun proved the following theorems \ref{Kane-Sun}-\ref{Kane-Sun(2)}, via modular forms and the theory of
quadratic forms [17]. Note that every positive integer $n$ can be expressed in the form $n=2^{v_{2}(n)}n'$ with
$v_{2}(n)\in\mathbb{N}$ and $n'$ odd. $v_{2}(a)$ is called the $2$-adic order of $a$ (equivalently
$2^{v_{2}(a)}\|a$) while $a'$ is said to be the odd part of $a$.\par\bigskip

\begin{teor}\label{Kane-Sun}
Fix $a,b,c\in\mathbb{Z}^{+}$ with $\mathrm{gcd}(a,b,c)=1$. Then the form \par\bigskip

\begin{centering}
$f(x,y,z)=at_{x}+bt_{y}+ct_{z}$\par\bigskip

\end{centering}

is asymptotically universal if and only if\par\bigskip

\begin{centering}

$-bc\hspace{0.1cm}R\hspace{0.1cm}a'$,\quad$-ac\hspace{0.1cm}R\hspace{0.1cm}b'$, and
$-ab\hspace{0.1cm}R\hspace{0.1cm}c'$.

\end{centering}

\end{teor}

\begin{teor}\label{Kane-Sun(1)}
Fix $a,b,c\in\mathbb{Z}^{+}$ with $\mathrm{gcd}(a,b,c)=1$. Then the form \par\bigskip

\begin{centering}
$f(x,y,z)=as_{x}+bt_{y}+ct_{z}$\par\bigskip

\end{centering}

is asymptotically universal if and only if we have the following (1)-(2)

\begin{enumerate}[(1)]
\item $-bc\hspace{0.1cm}R\hspace{0.1cm}a'$,\quad$-2ac\hspace{0.1cm}R\hspace{0.1cm}b'$,
and $-2ab\hspace{0.1cm}R\hspace{0.1cm}c'$.

\item Either $4\nmid b$ or $4\nmid c$.

\end{enumerate}

\end{teor}

\begin{teor}\label{Kane-Sun(2)}
Fix $a,b,c\in\mathbb{Z}^{+}$ with $\mathrm{gcd}(a,b,c)=1$. Then the form \par\bigskip

\begin{centering}
$f(x,y,z)=as_{x}+bs_{y}+ct_{z}$\par\bigskip

\end{centering}

is asymptotically universal if and only if we have the following (1)-(2)

\begin{enumerate}[(1)]
\item $-2bc\hspace{0.1cm}R\hspace{0.1cm}a'$,\quad$-2ac\hspace{0.1cm}R\hspace{0.1cm}b'$,
and $-ab\hspace{0.1cm}R\hspace{0.1cm}c'$.

\item Either $4\nmid c$, or both $4\| c$ and $2\|ab$.

\end{enumerate}

\end{teor}

Where if $E(f)=\{n\in\mathbb{N}\mid
f(x,y,z)=n\hspace{0.1cm}\mathrm{has\hspace{0.1cm}no\hspace{0.1cm}integral\hspace{0.1cm}solutions}\}$ has
asymptotic density zero then $f$ is \textit{asymptotically universal}, if $E(f)$ is finite then $f$ is
\textit{almost universal}. If $E(f)=\varnothing$, then $f$ is said to be \textit{universal}.\par\bigskip

$a\hspace{0.1cm}R\hspace{0.1cm}m$ if and only if the Legendre symbol $(\frac{a}{p})$ equals $1$ for every prime
divisor $p$ of $m$. That is, $a$ is quadratic residue modulo $m$. \par\bigskip

\setcounter{Nota}{7}

\begin{Nota}\label{forms}
For example each of the following forms represents every positive integer

\begin{enumerate}[(a)]
\item $f_{1}(u,v,w,x,y,z)$= $\alpha(t_{u}+4t_{v}+\beta(s_{w}+s_{w+1}))+(1-\alpha)(t_{x}+2s_{y}+2s_{z})$, $\alpha,\beta\in\{0,1\}$,\hspace{0.1cm} $u,x\geq0$,\hspace{0.1cm} $v\geq w\geq0$, $y\geq z\geq
1$,
\item $f_{2}(x,y,z)=t_{x}+t_{y}+2s_{z}$, $x,y,z\geq0$,

\item $f_{3}(x,y,z)=t_{x}+t_{y}+t_{z}$, $x,y,z\geq0$,

\item $f_{4}(x,y,z)=t_{x}+2s_{y}+4t_{z}$, $x,y,z\geq0$,

\item $f_{5}(x,y,z)=4t_{x}+t_{y}+t_{z}$, $x,y,z\geq0$.

\end{enumerate}

In fact Kane and Sun gave the complete lists of those forms $as_{x}+bs_{y}+ct_{z}$, $as_{x}+bt_{y}+ct_{z}$, with
$(a,b,c)\in\mathbb{Z}^{+}$ and $a+b+c\leq 10$ which are almost universal but not universal. In this case those
asymptotically universal ones are all almost universal.\par\bigskip The corresponding almost universal forms which
are not universal are respectively\par\bigskip

\begin{centering}
$\begin{array}{ccc}
  s_{x}+2s_{y}+3t_{z}, & 2s_{x}+4s_{y}+t_{z}, & s_{x}+6s_{y}+t_{z}, \\
  s_{x}+s_{y}+5t_{z},& 2s_{x}+3s_{y}+2t_{z}, & 3s_{x}+4s_{y}+t_{z}, \\
  s_{x}+2s_{y}+6t_{z}, & s_{x}+5s_{y}+3t_{z}, & 2s_{x}+4s_{y}+3t_{z}, \\
  4s_{x}+4s_{y}+t_{z}, & s_{x}+4s_{y}+5t_{z}, &
\end{array}$\par\bigskip

\end{centering}

\begin{centering}
$\begin{array}{ccc}
  5s_{x}+t_{y}+t_{z}\sim s_{x}+5s_{y}+2t_{z}, & 5s_{x}+2t_{y}+2t_{z}\sim 2s_{x}+5s_{y}+4t_{z}, & s_{x}+4t_{y}+2t_{z}, \\
  8s_{x}+t_{y}+t_{z}\sim s_{x}+8s_{y}+2t_{z}, & 2s_{x}+3t_{y}+2t_{z}, & 3s_{x}+4t_{y}+2t_{z}, \\
  2s_{x}+5t_{y}+t_{z}, &3s_{x}+5t_{y}+t_{z},  &5s_{x}+4t_{y}+t_{z},  \\
  4s_{x}+4t_{y}+t_{z}, &5s_{x}+3t_{y}+2t_{z},  &
\end{array}$\par\bigskip
\end{centering}

For the forms $at_{x}+bt_{y}+ct_{z}$ with $(a,b,c)\in\mathbb{Z}^{+}$ and $a+b+c\leq 10$, the following is the
complete list of those asymptotically universal forms which are not universal.\par\bigskip

 \begin{centering}
$\begin{array}{ccc}
  t_{x}+4t_{y}+4t_{z}\sim 4s_{x}+8t_{y}+t_{z}, &2t_{x}+3t_{y}+4t_{z},  &t_{x}+4t_{y}+5t_{z}, \\
  t_{x}+t_{y}+8t_{z}\sim s_{x}+8t_{y}+2t_{z}, & &  \\
  2t_{x}+2t_{y}+5t_{z}\sim2s_{x}+4t_{y}+5t_{z}, &  &  \\
  t_{x}+2t_{y}+6t_{z}. &  &
\end{array}$\par\bigskip
\end{centering}

Kane and sun conjectured that \par\bigskip

$E(s_{x}+2s_{y}+3t_{z})=\{23\}$,\quad $E(2s_{x}+4s_{y}+t_{z})=\{20\}$,\quad $E(s_{x}+5s_{y}+2t_{z})=\{19\}$,\quad
$E(s_{x}+6s_{y}+t_{z})=\{47\}$,\quad
$E(s_{x}+s_{y}+5t_{z})=\{3,11,12,27,129,138,273\}$,\quad$E(2s_{x}+3s_{y}+2t_{z})=\{1,19,43,94\}$,\quad
$E(2s_{x}+5s_{y}+t_{z})=\{4,27\}$,\quad $E(3s_{x}+4s_{y}+t_{z})=\{2,11,23,50,116,135,138\}$,\quad
$E(s_{x}+2s_{y}+6t_{z})=\{5,13,46,161\}$,\par\bigskip

$E(8s_{x}+t_{y}+t_{z})= E(s_{x}+8s_{y}+2t_{z})=\{5,40,217\}$,\quad $E(2s_{x}+3t_{y}+2t_{z})=\{1,16\}$,\quad
$E(2s_{x}+5t_{y}+t_{z})=\{4\}$,\quad $E(4s_{x}+3t_{y}+t_{z})=\{2,11,27,38,86,93,188,323\}$,\quad
$E(3s_{x}+5t_{y}+t_{z})=\{2,7\}$,\quad $E(3s_{x}+4t_{y}+2t_{z})=\{1,8,11,25\}$,\quad
$E(4s_{x}+4t_{y}+t_{z})=\{2,108\}$,\quad $E(6s_{x}+2t_{y}+t_{z})=\{4\}$,\quad
$E(5s_{x}+4t_{y}+t_{z})=\{2,16,31\}$,

\par\smallskip
 $E(5s_{x}+3t_{y}+2t_{z})=\{1,4,13,19,27,46,73,97,111,123,151,168\}$,
\par\bigskip

$E(2t_{x}+2t_{y}+5t_{z})=E(2s_{x}+4t_{y}+5t_{z})=\{1,3,10,16,28,43,46,85,169,175,211,223\}$,\par\bigskip
and\par\bigskip

$E(t_{x}+2t_{y}+6t_{z})=\{4,50\}$,\quad $E(2t_{x}+3t_{y}+4t_{z})=\{1,8,31\}$,\quad $E(t_{x}+4t_{y}+5t_{z})=\{2\}$.

\end{Nota}

\textbf{The main theorem}\par\bigskip

The following formulas (obtained by recursion) are solutions to the Guy's problem (see page 2) :\par\bigskip

For $n=8m+1$, $m\in\mathbb{N}$, we have that
\[n=
\begin{cases}
s_{2x+1}+4s_{2y+1}+4s_{2z+1},  &   \mathrm{if}\hspace{0.1cm}m=t_{x}+4t_{y}+s_{z}+s_{z+1},\hspace{0.1cm} x\geq0,\hspace{0.1cm} y\geq z\geq0. \\
s_{2x+1}+s_{4y}+s_{4z}, &\mathrm{if}\hspace{0.1cm}m=t_{x}+2s_{y}+2s_{z},\hspace{0.1cm}x\geq0,\hspace{0.1cm}y\geq
z\geq1.\\
s_{2x+1}+s_{4y}, &\mathrm{if}\hspace{0.1cm}m=t_{x}+2s_{y},\hspace{0.1cm}x\geq0,\hspace{0.1cm}y\geq0.
\end{cases}\]

If $n=8m+2$ then,

\[n=
\begin{cases}
s_{2x+1}+s_{2y+1},  &   \mathrm{if}\hspace{0.1cm}m=t_{x}+t_{y},\hspace{0.1cm} 0\leq x\leq y. \\
s_{4x}+s_{2y+1}+1, &\mathrm{if}\hspace{0.1cm}m=2s_{x}+t_{y},\hspace{0.1cm}x\geq1,\hspace{0.1cm}y\geq 0, \\
s_{2x+1}+s_{2y+1}+s_{4z} &\mathrm{if}\hspace{0.1cm}m=t_{x}+t_{y}+2s_{z},\hspace{0.1cm}1\leq x\leq
y,\hspace{0.1cm}z\geq1.
\end{cases}\]
\par\bigskip
For $n=8m+3$ we have that $n=s_{2x+1}+s_{2y+1}+s_{2z+1}$, if $m=t_{x}+t_{y}+t_{z}$, $x,y,z\geq0$.\par\bigskip If
$n=8m+5$ then,
\[n=
\begin{cases}
s_{2x+1}+s_{4y}+4s_{2z+1},  &   \mathrm{if}\hspace{0.1cm}m=t_{x}+2s_{y}+4t_{z},\hspace{0.1cm} x,z\geq0,\hspace{0.1cm} y\geq z+1. \\
s_{2x+1}+s_{4y+2}+s_{4z}, &\mathrm{if}\hspace{0.1cm}m=t_{x}+4t_{y}+2s_{z},\hspace{0.1cm}x\geq0,\hspace{0.1cm}y\geq
z\geq1.\\
s_{2x+1}+4s_{2z+1},  &   \mathrm{if}\hspace{0.1cm}m=t_{x}+4t_{z},\hspace{0.1cm} x,z\geq0.
\end{cases}\]

If $n=8m+6$ then,
\[n=
\begin{cases}
s_{4x+2}+s_{2y+1}+s_{2z+1},  &   \mathrm{if}\hspace{0.1cm}m=4t_{x}+t_{y}+t_{z},\hspace{0.1cm}x\geq0,\hspace{0.1cm} z\geq1,\hspace{0.1cm} y\geq z\geq1. \\
s_{4x+2}+s_{2y+1}+1, &\mathrm{if}\hspace{0.1cm}m=4t_{x}+t_{y},\hspace{0.1cm}x,y\geq0,\\
s_{2x+1}+s_{2y+1}+4, &\mathrm{if}\hspace{0.1cm}m=t_{x}+t_{y},\hspace{0.1cm}1\leq x\leq y.
\end{cases}\]

\par\bigskip

Since it is easy to verify that every number of the form $4^{k}(8m+7)$, $k,m\in\mathbb{N}$ cannot be expressible
as a sum of three or fewer square numbers [8], and if $n=4^{a}n_{1}$, $4\nmid n_{1}$ and $n_{1}$ is the sum of
three squares, say $n_{1}=\underset{n=1}{\overset{3}{\sum}}x^{2}_{i}$, then
$n=\underset{n=1}{\overset{3}{\sum}}x^{2}_{i}$ is also a sum of three squares. The formulas given above and
theorems \ref{Kane-Sun}-\ref{Kane-Sun(2)} provide a proof of the following result : \setcounter{teor}{8}
\begin{teor}\label{new-L-G}

If $t\notin\{n\in\mathbb{N}\hspace{0.1cm}\backslash\hspace{0.1cm}\{0\}\mid
n=4^{s}(8m+7),\hspace{0.1cm}\mathrm{for\hspace{0.1cm}some}\hspace{0.1cm}m,s\in\mathbb{N}\}$ then $t$ is the sum of
three or fewer squares.\hspace{0.5cm}\text{\qed}

\end{teor}

\end{document}